\documentclass[12pt]{article}                                               

\input{amssym.def}                                                           
\input{amssym.tex}                                                           
                                                                             
\begin{document}                                                             
\title{On the AF-algebra of a Hecke eigenform}

\author{Igor ~V. ~Nikolaev
\footnote{Partially supported 
by NSERC.}}


\date{} 
\maketitle


\newtheorem{thm}{Theorem}
\newtheorem{lem}{Lemma}
\newtheorem{dfn}{Definition}
\newtheorem{rmk}{Remark}
\newtheorem{cor}{Corollary}
\newtheorem{prp}{Proposition}
\newtheorem{exm}{Example}
\newtheorem{cnj}{Conjecture}

\newcommand{\N}{{\Bbb N}}
\newcommand{\mod}{\hbox{\bf mod}}
\newcommand{\GCD}{\hbox{\bf GCD}}
\newcommand{\GL}{\hbox{\bf GL}}
\newcommand{\rank}{\hbox{\bf rank}}
\newcommand{\End}{\hbox{\bf End}}
\newcommand{\Per}{\hbox{\bf Per}}
\newcommand{\ch}{\hbox{\bf char}}

\begin{abstract}
An AF-algebra is assigned to each cusp form $f$ of weight two;
we study properties of this operator algebra, when $f$ is a Hecke  eigenform.

\vspace{7mm}


\vspace{5mm}
{\it AMS (MOS) Subj. Class.: 11F03;  46L85}
\end{abstract}

\section{Introduction}
\noindent
{\bf A. The Modularity Theorem} asserts  that:

\medskip\hskip2cm
{\it All rational elliptic curves arise from modular forms.}

\medskip\noindent
This  result  is tremendously important,  since it leads to a spectacular
proof of  Fermat's Last Theorem.   The reader  of the excellent book \cite{DS} on the Modularity
Theorem  will find on page (x) of the introduction the following interesting object.  
(We shall modify  the original text to match  our notation, which can be found
in Section 2.)  
Denote by $f_N\in S_2(\Gamma_0(N))$ a Hecke eigenform and by $f_N^{\sigma}$ all its
conjugates;  consider a lattice $\Lambda_{f_N}$ generated by the complex periods
of holomorphic forms $\omega_N^{\sigma}=f_N^{\sigma}dz$ on the Riemann surface 
 $X_0(N)={\Bbb H}^*/\Gamma_0(N)$.  If $|\sigma|$ is the number of conjugates, 
the abelian variety  ${\cal A}_{f_N}:={\Bbb C}^{|\sigma|}/\Lambda_{f_N}$ is said to be {\it associated}
to the eigenform $f_N$; it has the following remarkable property (the Modularity
Theorem):

\medskip\hskip0.5cm
{\it There exists a homomorphism of ${\cal A}_{f_N}$ onto a rational elliptic curve.}

\medskip\noindent 
Let $\phi_N=\Re~(\omega_N)$ be the real part
of $\omega_N$; it is a closed form on the surface $X_0(N)$.
(Alternatively, one can take for $\phi_N$ the  imaginary part of $\omega_N$.)
Clearly, $\omega_N$ defines a unique form $\phi_N$;  the converse follows from the
Hubbard-Masur Theorem \cite{HuMa1}.  Since $\omega_N$ and $\phi_N$
define each other,  what object will replace the associated variety ${\cal A}_{f_N}$
in the case of $\phi_N$?    
Roughly speaking, it is shown in this paper  that such a replacement is given 
by an operator algebra ${\goth A}_{f_N}$ coming  from  the  real  periods of the form $\phi_N$; 
we study the basic properties of such an algebra (Theorem \ref{thm1}).

\medskip\noindent
{\bf B. The AF-algebra  ${\goth A}_{f_N}$.}
Let $f\in S_2(\Gamma_0(N))$ be a cusp form and $\omega=fdz$
the corresponding holomorphic differential on $X_0(N)$.  We shall 
denote by $\phi=\Re~(\omega)$  a  closed  form on $X_0(N)$
and  consider its periods $\lambda_i=\int_{\gamma_i}\phi$
against  a basis $\gamma_i$  in  the (relative) homology group 
$H_1(X_0(N), Z(\phi); ~{\Bbb Z})$,  where $Z(\phi)$ is the set of zeros of $\phi$.   
Assume $\lambda_i> 0$ and consider the vector $\theta=(\theta_1,\dots,\theta_{n-1})$
with $\theta_i=\lambda_{i+1} / \lambda_1$. The Jacobi-Perron continued fraction of
$\theta$ (\cite{BE}) is given by the formula:
$$
\left(\matrix{1\cr \theta}\right)=
\lim_{i\to\infty} \left(\matrix{0 & 1\cr I & b_1}\right)\dots
\left(\matrix{0 & 1\cr I & b_i}\right)
\left(\matrix{0\cr {\Bbb I}}\right)=
\lim_{i\to\infty} B_i\left(\matrix{0\cr {\Bbb I}}\right),
$$
where $b_i=(b^{(i)}_1,\dots, b^{(i)}_{n-1})^T$ is a vector of  non-negative integers,  
$I$ is the unit matrix and ${\Bbb I}=(0,\dots, 0, 1)^T$.
By ${\goth A}_{f}$ we shall understand  the  AF-algebra
given its  Bratteli diagram  with  partial multiplicity matrices $B_i$. 
Recall  that an  AF-algebra is called  {\it stationary}  if $B_i=B=$ Const \cite{E}.
When two non-similar matrices $B$ and $B'$ have the same characteristic polynomial, 
the corresponding stationary AF-algebras will be called   {\it companion AF-algebras}.   
Denote by ${\goth A}_{f_N}$  an AF-algebra, such that  $f_N\in S_2(\Gamma_0(N))$
is a Hecke eigenform. Our main result can be stated  as follows.
\begin{thm}\label{thm1}
The AF-algebra  ${\goth A}_{f_N}$ is stationary  unless $f_N$ is a rational eigenform,
in which case ${\goth A}_{f_N}\cong {\Bbb C}$; moreover,  ${\goth A}_{f_N}$ and  ${\goth A}_{f_N^{\sigma}}$
are companion AF-algebras.  
\end{thm}
The paper is organized as follows. The minimal preliminary results are
expounded in Section 2, where  we introduce the Hecke eigenforms, the AF-algebras
and the Jacobi-Perron continued fractions.  Theorem \ref{thm1} is proved in Section 3.

\section{Preliminaries}
{\bf A. The Hecke eigenforms.}
Let $N>1$ be a natural number and  consider a (finite index) subgroup 
of the modular group given by the formula: 
$$
\Gamma_0(N) = \left\{\left(\matrix{a & b\cr c & d}\right)\in SL(2,{\Bbb Z})~|~
c\equiv 0~\mod ~N\right\}.
$$
Let ${\Bbb H}=\{z=x+iy\in {\Bbb C} ~|~ y>0\}$ be the upper half-plane  and 
let $\Gamma_0(N)$  act on ${\Bbb H}$  by the linear fractional
transformations;  consider an orbifold  ${\Bbb H}/\Gamma_0(N)$.
To compactify the orbifold 
at the cusps, one adds a boundary to ${\Bbb H}$,  so that 
${\Bbb H}^*={\Bbb H}\cup {\Bbb Q}\cup\{\infty\}$ and the compact Riemann surface 
$X_0(N)={\Bbb H}^*/\Gamma_0(N)$ is called a {\it modular curve}.   
The meromorphic functions $f(z)$ on ${\Bbb H}$ that
vanish at the cusps and such that
$$
f\left({az+b\over cz+d}\right)={1\over (cz+d)^2}f(z),\qquad
\forall \left(\matrix{a & b\cr c & d}\right)\in\Gamma_0(N), 
$$
are  called  {\it cusp forms} of weight two;  the (complex linear) space of such forms
will be denoted by $S_2(\Gamma_0(N))$.  The formula $f(z)\mapsto \omega=f(z)dz$ 
defines an isomorphism  $S_2(\Gamma_0(N))\cong \Omega_{hol}(X_0(N))$, where 
$\Omega_{hol}(X_0(N))$ is the space of holomorphic differentials
on the Riemann surface $X_0(N)$.  Note that 
\linebreak
$\dim_{\Bbb C}(S_2(\Gamma_0(N))=\dim_{\Bbb C}(\Omega_{hol}(X_0(N))=g$,
where $g=g(N)$ is the genus of the surface $X_0(N)$. 
A Hecke operator, $T_n$, acts on $S_2(\Gamma_0(N))$ by the formula
$T_n f=\sum_{m\in {\Bbb Z}}\gamma(m)q^m$, where
$\gamma(m)= \sum_{a|\GCD(m,n)}a c_{mn/a^2}$ and 
$f(z)=\sum_{m\in {\Bbb Z}}c(m)q^m$ is the Fourier
series of the cusp form $f$ at $q=e^{2\pi iz}$.  Further,  $T_n$ is a
self-adjoint linear operator on the vector space $S_2(\Gamma_0(N))$
endowed with the Petersson inner product;  the algebra
${\Bbb T}_N :={\Bbb Z}[T_1,T_2,\dots]$ is a commutative algebra.
Any cusp form $f_N\in S_2(\Gamma_0(N))$ that is an eigenvector
for one (and hence all) of $T_n$, is referred  to
as a {\it Hecke eigenform}; such an  eigenform is called
{\it rational}  whenever its Fourier coefficients $c(m)\in {\Bbb Z}$.
The Fourier  coefficients $c(m)$ of $f_N$ are algebraic integers,  and we 
denote by ${\Bbb K}_{f_N}={\Bbb Q}(c(m))$ an extension of the field ${\Bbb Q}$ 
by the Fourier coefficients of $f_N$. Then ${\Bbb K}_{f_N}$
is a real algebraic number field of degree $1\le \deg~({\Bbb K}_{f_N} / {\Bbb Q})\le g$,
where $g$ is the genus of the surface $X_0(N)$ \cite{DS}, Proposition 6.6.4. 
Any embedding $\sigma: {\Bbb K}_{f_N}\to {\Bbb C}$ conjugates $f_N$ by acting
on its coefficients;  we write the corresponding Hecke  eigenform  
$f_N^{\sigma}(z):=\sum_{m\in {\Bbb Z}}\sigma(c(m))q^m$.

\medskip\noindent
{\bf B. The AF-algebras.}
A {\it $C^*$-algebra} is an algebra $A$ over ${\Bbb C}$ with a norm
$a\mapsto ||a||$ and an involution $a\mapsto a^*$ such that
it is complete with respect to the norm and $||ab||\le ||a||~ ||b||$
and $||a^*a||=||a^2||$ for all $a,b\in A$.
Any commutative $C^*$-algebra is  isomorphic
to the algebra $C_0(X)$ of continuous complex-valued
functions on some locally compact Hausdorff space $X$; 
otherwise, $A$ represents a noncommutative  topological
space.  
The $C^*$-algebras $A$ and $A'$ are said to be
{\it stably isomorphic} (Morita equivalent)  if   
$A\otimes {\cal K}\cong A'\otimes {\cal K}$, where ${\cal K}$
is the $C^*$-algebra of compact operators; 
roughly speaking,  stable isomorphism 
means that $A$ and $A'$ are
homeomorphic as noncommutative topological spaces.

An {\it AF-algebra}  (Approximately Finite $C^*$-algebra) is defined to
be the  norm closure of an ascending sequence of   finite dimensional
$C^*$-algebras $M_n$,  where  $M_n$ is the $C^*$-algebra of the $n\times n$ matrices
with entries in ${\Bbb C}$. Here the index $n=(n_1,\dots,n_k)$ represents
the  semi-simple matrix algebra $M_n=M_{n_1}\oplus\dots\oplus M_{n_k}$.
The ascending sequence mentioned above  can be written as 
$M_1\buildrel\rm\varphi_1\over\longrightarrow M_2
   \buildrel\rm\varphi_2\over\longrightarrow\dots,
$
where $M_i$ are the finite dimensional $C^*$-algebras and
$\varphi_i$ the homomorphisms between such algebras.  
The homomorphisms $\varphi_i$ can be arranged into  a graph as follows. 
Let  $M_i=M_{i_1}\oplus\dots\oplus M_{i_k}$ and 
$M_{i'}=M_{i_1'}\oplus\dots\oplus M_{i_k'}$ be 
the semi-simple $C^*$-algebras and $\varphi_i: M_i\to M_{i'}$ the  homomorphism. 
One has  two sets of vertices $V_{i_1},\dots, V_{i_k}$ and $V_{i_1'},\dots, V_{i_k'}$
joined by  $b_{rs}$ edges  whenever the summand $M_{i_r}$ contains $b_{rs}$
copies of the summand $M_{i_s'}$ under the embedding $\varphi_i$. 
As $i$ varies, one obtains an infinite graph called the  {\it Bratteli diagram} of the
AF-algebra.  The matrix $B=(b_{rs})$ is known as  a {\it partial multiplicity matrix};
an infinite sequence of $B_i$ defines a unique AF-algebra.

For a unital $C^*$-algebra $A$, let $V(A)$
be the union (over $n$) of projections in the $n\times n$
matrix $C^*$-algebra with entries in $A$;
projections $p,q\in V(A)$ are {\it equivalent} if there exists a partial
isometry $u$ such that $p=u^*u$ and $q=uu^*$. The equivalence
class of projection $p$ is denoted by $[p]$;
the equivalence classes of orthogonal projections can be made to
a semigroup by putting $[p]+[q]=[p+q]$. The Grothendieck
completion of this semigroup to an abelian group is called
the  $K_0$-group of the algebra $A$.
The functor $A\to K_0(A)$ maps the category of unital
$C^*$-algebras into the category of abelian groups, so that
projections in the algebra $A$ correspond to a positive
cone  $K_0^+\subset K_0(A)$ and the unit element $1\in A$
corresponds to an order unit $u\in K_0(A)$.
The ordered abelian group $(K_0,K_0^+,u)$ with an order
unit  is called a {\it dimension group}; an order-isomorphism
class of the latter we denote by $(G,G^+)$.

\medskip\noindent
{\bf C. The Jacobi-Perron fractions.}
Let $a_1,a_2\in {\Bbb N}$ such that $a_2\le a_1$. Recall that the greatest common
divisor of $a_1,a_2$, $\GCD(a_1,a_2)$, can be determined from the Euclidean algorithm:
$$
\left\{
\begin{array}{cc}
a_1 &= a_2b_1 +r_3\nonumber\\
a_2 &= r_3b_2 +r_4\nonumber\\
r_3 &= r_4b_3 +r_5\nonumber\\
\vdots & \nonumber\\
r_{k-3} &= r_{k-2}b_{k-1}+r_{k-1}\nonumber\\
r_{k-2} &= r_{k-1}b_k,
\end{array}
\right.
$$
where $b_i\in {\Bbb N}$ and $\GCD(a_1,a_2)=r_{k-1}$. 
The Euclidean algorithm can be written as the regular continued 
fraction
$$
\theta={a_1\over a_2}=b_1+{1\over\displaystyle b_2+
{\strut 1\over\displaystyle +\dots+ {1\over b_k}}}
=[b_1,\dots b_k].
$$
If $a_1, a_2$ are non-commensurable  in the sense that $\theta\in {\Bbb R}-{\Bbb Q}$,
then the Euclidean algorithm never stops,  and $\theta=[b_1, b_2, \dots]$. Note that the regular  
continued fraction can be written in  matrix form
$$
\left(\matrix{1\cr \theta}\right)=
\lim_{k\to\infty} \left(\matrix{0 & 1\cr 1 & b_1}\right)\dots
\left(\matrix{0 & 1\cr 1 & b_k}\right)
\left(\matrix{0\cr 1}\right). 
$$
The Jacobi-Perron algorithm and connected (multidimensional) continued 
fraction generalizes the Euclidean algorithm to the case $\GCD(a_1,\dots,a_n)$
when $n\ge 2$. Namely, let $\lambda=(\lambda_1,\dots,\lambda_n)$,
$\lambda_i\in {\Bbb R}-{\Bbb Q}$ and  $\theta_{i-1}={\lambda_i\over\lambda_1}$, where
$1\le i\le n$.   The continued fraction 
$$
\left(\matrix{1\cr \theta_1\cr\vdots\cr\theta_{n-1}} \right)=
\lim_{k\to\infty} 
\left(\matrix{0 &  0 & \dots & 0 & 1\cr
              1 &  0 & \dots & 0 & b_1^{(1)}\cr
              \vdots &\vdots & &\vdots &\vdots\cr
              0 &  0 & \dots & 1 & b_{n-1}^{(1)}}\right)
\dots 
\left(\matrix{0 &  0 & \dots & 0 & 1\cr
              1 &  0 & \dots & 0 & b_1^{(k)}\cr
              \vdots &\vdots & &\vdots &\vdots\cr
              0 &  0 & \dots & 1 & b_{n-1}^{(k)}}\right)
\left(\matrix{0\cr 0\cr\vdots\cr 1} \right),
$$
where $b_i^{(j)}\in {\Bbb N}\cup\{0\}$, is called the {\it Jacobi-Perron
algorithm (JPA)}. Unlike the regular continued fraction algorithm,
the JPA may diverge for certain vectors $\lambda\in {\Bbb R}^n$. However, 
for points of a generic subset of ${\Bbb R}^n$, the JPA converges \cite{Bau1}; 
in particular, the JPA for periodic fractions is always convergent.

\section{Proof of theorem \ref{thm1}}
 A standard dictionary (\cite{E}) between AF-algebras and their
dimension groups is adopted.  Instead of dealing with ${\goth A}_f$, we work with its dimension group 
$G_{{\goth A}_f}=(G, G^+)$, where $G\cong {\Bbb Z}^n$ is the lattice and 
$G^+=\{(x_1,\dots, x_n)\in {\Bbb Z}^n~|~\theta_1x_1+\dots+\theta_{n-1} x_{n-1}+ x_n\ge 0\}$ 
is a positive cone.  Recall, that  $G_{{\goth A}_f}$ is  abelian group with an order, which  
defines the AF-algebra ${\goth A}_f$,  up to a stable isomorphism.   
We arrange the proof in a series of lemmas. First, let us show,  that
${\goth A}_f$  is a correctly-defined AF-algebra.
\begin{lem}\label{lm1}
The ${\goth A}_f$ does not depend, up to a stable isomorphism,  on 
a basis in $H_1(X_0(N), Z(\phi); ~{\Bbb Z})$.
\end{lem}
{\it Proof.} 
Denote by ${\goth m}:= {\Bbb Z}\lambda_1+\dots+{\Bbb Z}\lambda_n$
a ${\Bbb Z}$-module in the real line ${\Bbb R}$. Let $\{\gamma_i'\}$ be a new basis
in $H_1(X_0(N), Z(\phi); {\Bbb Z})$, such that $\gamma_i'=\sum_{j=1}^na_{ij}\gamma_j$
for  matrix $A=(a_{ij})\in \GL_n({\Bbb Z})$. Using the integration rules, one gets:
$\lambda_i'  = \int_{\gamma_i'}\phi = \int_{\sum_{j=1}^na_{ij}\gamma_j}\phi=
 \sum_{j=1}^n\int_{\gamma_j}\phi  = \sum_{j=1}^na_{ij}\lambda_j$.
Thus, ${\goth m}'={\goth m}$ and a  change of basis in the homology group
$H_1(X_0(N), Z(\phi); {\Bbb Z})$ amounts to a change of basis in the module
${\goth m}$. It is an easy exercise to show  that there exists a linear transformation of ${\Bbb Z}^n$
sending the positive cone $G^+$ of $G_{{\goth A}_f}$ to the positive cone $(G^{+})'$ of 
$G_{{\goth A}_f'}$.  In other words, ${\goth A}_f'$ and  ${\goth A}_f$ are stably isomorphic.
$\square$

\begin{lem}\label{lm2}
The (scaled) periods $\lambda_i$ belong to the field ${\Bbb K}_{f_N}$. 
\end{lem}
{\it Proof.}  
Let ${\goth m}={\Bbb Z}\lambda_1+\dots+{\Bbb Z}\lambda_{2g}$ be a
${\Bbb Z}$-module generated by $\lambda_i$;  we seek  the
effect of the Hecke operators $T_m$ on ${\goth m}$.  By the
definition of a Hecke eigenform, $T_mf_N=c(m)f_N$ for all $T_m\in {\Bbb T}_N$.
In view of the isomorphism $S_2(\Gamma_0(N))\cong \Omega_{hol}(X_0(N))$,
one gets $T_m\omega_N=c(m)\omega_N$, where  $\omega_N=f_N dz$. Then 
$\Re~(T_m\omega_N)=T_m(\Re~(\omega_N))=\Re~(c(m)\omega_N)=c(m)  \Re~(\omega_N)$.
Therefore, $T_m\phi_N=c(m)\phi_N$, where $\phi_N=\Re~(\omega_N)$. The action of $T_m$
on ${\Bbb Z}$-module ${\goth m}$ can be written as 
$T_m({\goth m})=\int_{H_1}T_m\phi_N=\int_{H_1}c(m)\phi_N=c(m) {\goth m}$,
where $H_1:= H_1(X_0(N), Z(\phi_N); {\Bbb Z})$.  Thus, the Hecke operator $T_m$
acts  on the module ${\goth m}$ as multiplication by an algebraic
integer $c(m)\in {\Bbb K}_{f_N}$.

The action of $T_m$ on ${\goth m}={\Bbb Z}\lambda_1+\dots+{\Bbb Z}\lambda_n$
can be written as $T_m\lambda=c(m)\lambda$, where $\lambda=(\lambda_1,\dots,\lambda_n)$;
thus, $T_m$ is a linear operator (on the space ${\Bbb R}^n$), whose eigenvector $\lambda$
corresponds to the eigenvalue $c(m)$. It is an easy exercise in linear algebra  that 
$\lambda$ can be scaled so that all $\lambda_i$ lie in the same field as $c(m)$;
lemma \ref{lm2} follows.   
$\square$

\bigskip
{\sf Case I.} 
Let $f_N$ be not a  rational eigenform; then 
$n=\deg~({\Bbb K}_{f_N}/{\Bbb Q})\ge 2$.  Note,  that 
${\goth m}={\Bbb Z}\lambda_1+\dots+{\Bbb Z}\lambda_n$ is a full (i.e. the maximal
rank) ${\Bbb Z}$-module in the number field ${\Bbb K}_{f_N}$. 
Indeed, $\rank ~({\goth m})$ cannot exceed $n$, since ${\goth m}\subset {\Bbb K}_{f_N}$
and ${\Bbb K}_{f_N}$ is a vector space (over ${\Bbb Q}$) of dimension $n$.
On the other hand, $(\lambda_1,\dots,\lambda_n)$ is a basis of the field ${\Bbb K}_{f_N}$
and, as such, $\rank~({\goth m})$ cannot be less than $n$; thus, $\rank~({\goth m})=n$.
\begin{lem}\label{lm3}
The vector $(\lambda_1,\dots,\lambda_n)$ has a  periodic (Jacobi-Perron) continued fraction.
\end{lem}
{\it Proof.}
Since ${\goth m}\subset {\Bbb K}_{f_N}$ is a full ${\Bbb Z}$-module,
its endomorphism ring $\End~({\goth m})=\{\alpha\in {\Bbb K}_{f_N} : \alpha {\goth m}\subseteq {\goth m}\}$
is an order (a subring of the ring of integers) of the number field ${\Bbb K}_{f_N}$;  let $u$ be a unit of 
 the order \cite{BS}, p 112.  The action of $u$ on ${\goth m}$ can be written in a matrix form
$A\lambda=u\lambda$, where $\lambda$ is a basis in ${\goth m}$ and $A\in\GL_n({\Bbb Z})$;
with no loss of generality,  one can assume the matrix $A$ to be non-negative in a proper basis
of ${\goth m}$.

According to \cite{Bau1}, Prop.3, the matrix $A$
can be uniquely factorized as 
$A=\left(\small\matrix{0 & 1\cr I & b_1}\right)\dots\left(\matrix{0 & 1\cr I & b_k}\right)$,
where vectors $b_i=(b^{(i)}_1,\dots, b^{(i)}_{n-1})^T$ have  non-negative integer entries.
By \cite{Per1}, Satz XII, the periodic continued fraction
\begin{equation}\label{eq1}
\left(\matrix{1\cr \theta'}\right)=
\Per
~\overline{
 \left(\matrix{0 & 1\cr I & b_1}\right)\dots
\left(\matrix{0 & 1\cr I & b_k}\right)
}
\left(\matrix{0\cr {\Bbb I}}\right)
\end{equation}
converges to a vector $\lambda'=(\lambda_1',\dots,\lambda_n')$
which satisfies the equation $A\lambda'=u\lambda'$.
Since $A\lambda=u\lambda$, the vectors $\lambda$ and $\lambda'$ are collinear; 
but collinear vectors have the same continued fractions \cite{BE}. 
$\square$

\bigskip
The first case of Theorem \ref{thm1} follows from lemma \ref{lm3},
since ${\goth A}_{f_N}$ is a stationary  AF-algebra, whose period is given by the matrix 
$A$.

\bigskip
{\sf Case II.} 
Let $f_N$ be a rational eigenform; in this case $n=1$ and ${\Bbb K}_{f_N}={\Bbb Q}$.
The Bratteli diagram of ${\goth A}_{f_N}$ is finite and one-dimensional;
therefore, ${\goth A}_{f_N}\cong M_1({\Bbb C})={\Bbb C}$. This argument finishes the proof of
the first part of Theorem \ref{thm1}.

\bigskip
To prove the second part, 
let on the contrary $A\ne A'$ be similar matrices. To find $S$ such that 
$A'=S^{-1}AS$,  notice that   ${\goth m}^{\sigma}=\lambda_1^{\sigma}{\Bbb Z}+\dots+\lambda_n^{\sigma}{\Bbb Z}$.
Since ${\goth m}^{\sigma}={\goth m}$,  $\lambda_j^{\sigma}=\sum s_{ij}\lambda_i$, where $S=(s_{ij})$;
but $\sigma^k=$Id  for some integer $k$ and  thus $S^k=I$. Therefore, $(A')^k=(S^{-1}AS)^k=A^k$
and $A'=A$, which contradicts our assumption. On the other hand, $\lambda_j^{\sigma}\in {\Bbb K}_{f_N}$
implies that the characteristic polynomials $\ch~(A)=\ch~(A')$; therefore, ${\goth A}_{f_N}$ and ${\goth A}_{f_N^{\sigma}}$
are companion AF-algebras. 
$\square$

\bigskip\noindent
{\sf Acknowledgments.}  
I am grateful to Prof. ~Yu.~I.~Manin for helpful remarks and to Prof. ~D.~Solomon for
his interest;  Dr. L.~D.~Taylor kindly provided me with a copy of \cite{Tay1}.  
I  thank the referee for prompt and thoughtful comments.  



\vskip1cm

\textsc{The Fields Institute for Mathematical Sciences, Toronto, ON, Canada,  
E-mail:} {\sf igor.v.nikolaev@gmail.com}

\smallskip
{\it Current address: 616-315 Holmwood Ave., Ottawa, ON, Canada, K1S 2R2}

\end{document}